\documentclass{amsart}

\usepackage{graphicx}
\usepackage{amssymb, comment, color}
\usepackage{tikz-cd}
\usepackage{authblk}

\newcommand{\N}{\mathbb{N}}                     
\newcommand{\Z}{\mathbb{Z}}                     
\newcommand{\R}{\mathbb{R}}                     
\newcommand{\C}{\mathbb{C}}                     
\newcommand{\U}{\mathcal{U}}                     
\newcommand{\F}{\mathcal{F}}                     

\renewcommand{\S}{\mathcal{S}}
\newcommand{\im}{\mathrm{Im\,}}                 
\newcommand{\re}{\mathrm{Re\,}}                 
\newcommand{\CP}{\mathbb{CP}}                   
\newcommand{\RP}{\mathbb{RP}}                   

\newcommand{\cz}{{\rm CZ}}
\newcommand{\wind}{{\rm wind}}

\newtheorem{thm}{\sc Theorem}[section]          
\newtheorem*{thm*}{\sc Theorem}               	
\newtheorem*{cor*}{\sc Corollary}        		
\newtheorem{prop}[thm]{\sc Proposition}     
\newtheorem{defn}[thm]{\sc Definition}      
\newtheorem{rem}[thm]{\sc Remark}           

\newtheorem{theo}[thm]{\sc Theorem}
\newtheorem{defi}[thm]{\sc Definition}

\usepackage{lipsum}

\usepackage[foot]{amsaddr}

\usepackage{amsfonts,amsthm,amssymb,latexsym,amsmath,setspace,amscd,verbatim}
\usepackage[matrix,arrow]{xy}
\usepackage[active]{srcltx}
\usepackage[colorlinks=true]{hyperref}
\usepackage{mathrsfs}

\usepackage{authblk}

\begin{document}

\title[Reeb flows and Transverse Foliations]{Reeb flows, pseudo-holomorphic curves and transverse foliations}

\author{Naiara V. de Paulo}

\address{Naiara V. de Paulo, Universidade Federal de Santa Catarina, Departamento de Matem\'atica, Rua Jo\~ao Pessoa, 2514, Bairro Velha, Blumenau SC, Brazil 89036-004}
\email{naiara.vergian@ufsc.br}

\author{Pedro A. S. Salom\~ao}

\address{Pedro A. S. Salom\~ao, Instituto de Matemática e Estatística, Departamento de Matemática, Universidade de São Paulo, Rua do Matão, 1010, Cidade Universitária, São
Paulo, SP 05508-090, Brazil, and NYU-ECNU Institute of Mathematical Sciences at NYU Shanghai, 3663 Zhongshan Road North, Shanghai, 200062, China}
\email{psalomao@ime.usp.br, pas383@nyu.edu}

\begin{abstract}
Pseudo-holomorphic curves in symplectizations, as introduced by Hofer in 1993, and then developed by Hofer, Wysocki, and Zehnder, have brought new insights to Hamiltonian dynamics, providing new approaches to some classical questions in Celestial Mechanics. This short survey presents some recent developments in Reeb dynamics based on the theory of pseudo-holomorphic curves in symplectizations, focusing on transverse foliations near critical energy surfaces.

\end{abstract}

\dedicatory
{On the occasion of the Golden Jubilee Celebration of IME-USP}

\maketitle


\section{Introduction}

The early twentieth century is a milestone in developing the qualitative theory of dynamical systems with significant contributions by Birkhoff and Poincar\'e to the study of geodesic flows on surfaces and the restricted three-body problem. Global surfaces of section, called Poincar\'e sections, were used to tackle questions concerning the existence of periodic motions, boosting the development of surface dynamics and culminating with the construction of a vast range of algebraic invariants. The Arnold conjecture, on the least number of fixed points of Hamiltonian diffeomorphisms, triggered the use of pseudo-holomorphic curves to build homology theories that detect periodic orbits and other dynamical objects such as transverse foliations, a notion that generalizes Poincar\'e sections.

Let $\phi_t,t\in \R,$ be a smooth flow on a smooth closed $3$-manifold $M$. An embedded compact surface $\Sigma \hookrightarrow M$ is called a global surface of section for $\phi_t$ if the boundary of $\Sigma$ is formed by finitely many periodic orbits, the flow is transverse to $\dot \Sigma = \Sigma \setminus \partial \Sigma$ and every trajectory in $M\setminus \partial \Sigma$ hits $\dot \Sigma$  in the past and in the future. The diffeomorphism  $\psi: \dot \Sigma \to \dot \Sigma$ given by the first return map encodes the dynamics of $\phi_t$. 

In \cite{Po}, Poincar\'e found annulus-like global surfaces of section in the circular planar restricted three-body problem for low energies and small mass ratios. He explored the twist condition on the first return map to obtain infinitely many periodic orbits. Eventually, Birkhoff established a more general fixed point theorem for area-preserving homeomorphisms of the annulus, which afterwards became known as the Poincar\'e-Birkhoff theorem.

In the 1960s, Arnold interpreted the Poincar\'e-Birkhoff theorem as a fixed point theorem for Hamiltonian diffeomorphisms on the $2$-torus. This interpretation led Arnold to conjecture that any Hamiltonian diffeomorphism of a closed symplectic manifold contains a minimum number of fixed points that depends only on the manifold's topology.

In the 1980s, Conley and Zehnder \cite{CZ} proved the Arnold conjecture for $2n$-dimensional tori equipped with the standard symplectic form. Floer \cite{floer} used pseudo-holomorphic curves -- introduced by Gromov \cite{gromov} in his proof of the non-squeezing theorem --  to confirm the Arnold conjecture in a broader range of symplectic manifolds.

In 1993 Helmut Hofer \cite{93} used finite-energy pseudo-holomorphic curves in symplectizations to study Hamiltonian dynamics restricted to contact-type energy surfaces, proving many instances of the Weinstein conjecture in dimension three.  Hofer, Wysocki, and Zehnder then paved the foundations of the theory of such curves in a series of articles \cite{props1,props2,props3,convex,fols}, obtaining some profound results in Hamiltonian dynamics:
\begin{itemize}
\item every strictly convex hypersurface in $\R^4$ admits a disk-like global surface of section bounded by a periodic orbit with Conley-Zehnder index $3$, see \cite{convex}.
\item every generic star-shaped hypersurface in $\R^4$ admits a transverse foliation  whose binding orbits have Conley-Zehnder indices $1$, $2$ or $3$, see \cite{fols}.
\end{itemize}

The theory of pseudo-holomorphic curves in symplectizations provided a new understanding of Reeb dynamics in dimension $3$ and, in particular, established many versions of contact homologies based on Floer's theoretical methods.

In this survey, we compile some recent developments in three-dimensional Reeb dynamics obtained with pseudo-holomorphic curves. We focus on the existence of  transverse foliations, particularly for Hamiltonian flows on energy surfaces in a neighborhood of a critical level. 

\section*{Acknowledgements} P. Salom\~ao acknowledges the support of NYU-ECNU Institute of Mathematical Sciences at NYU Shanghai. P. Salomão is partially supported by FAPESP 2016/25053-8 and CNPq 306106/2016-7.

\subsection{Basic notions in Reeb dynamics}\label{sec_basics}
Reeb flows compose a class of conservative dynamical systems that includes geodesic flows and many models in Celestial Mechanics  \cite{AFvKP}. 

A contact form  on a smooth oriented closed $3$-manifold $M$ is a  $1$-form $\lambda$ satisfying $\lambda \wedge d\lambda\neq 0.$ It determines a tangent plane distribution
$$
\xi=\ker \lambda \subset TM,
$$
called the contact structure, and a vector field $R$ on $M$, determined by 
\begin{equation}\label{eq_reeb}
 d \lambda(R, \cdot)=0 \ \mbox{ and } \ \lambda (R) = 1,
\end{equation}
 called the Reeb vector field of $\lambda$. 

If a hypersurface $i:S\hookrightarrow \R^4$ is star-shaped with respect to the origin, then $i^*\Lambda_0$ is a contact form on $S$, where $\Lambda_0$ is the Liouville form 
\begin{equation}\label{eq_liou}
\Lambda_0 = \frac{1}{2} \sum_{i=1}^2 (x_idy_i - y_i dx_i).
\end{equation}  The simplest example is the $3$-sphere $i_0:S^3 = \{x_1^2+y_1^2+x_2^2+y_2^2=1\}\hookrightarrow \R^4$, whose contact form $\lambda_0= i_0^*\Lambda_0$ is called standard and whose contact structure $\xi_0 = \ker  \lambda_0$ is called tight.

Hamiltonian dynamics restricted to contact-type energy levels are Reeb dynamics. Let $(W,\omega)$ be  a symplectic manifold and let $H: W \to \R$ be a smooth function. The Hamiltonian vector field $X_H$ is uniquely determined by  $\iota_{X_H} \omega = dH$. 
A regular energy level $i:M=H^{-1}(c)\hookrightarrow W$  has contact-type if there exists a Liouville vector field $Y$ (i.e. $\mathcal{L}_Y \omega=\omega$) defined in a neighborhood  of $M$ such that $Y$ is everywhere transverse to $M$. In this case, 
$\lambda:= i^*(\iota_Y \omega)$
is a contact form on $M$ satisfying $d\lambda=i^*\omega$. In this case, 
$X_H|_M \subset \ker (d\lambda)$ is  parallel to the Reeb vector field of $\lambda$.

The Reeb flow $\psi_t:M\to M,t\in \R,$  preserves the contact form $\lambda$. In particular, it preserves the contact structure $\xi=\ker \lambda$ and the fiberwise symplectic form $d\lambda|_\xi$.



A fundamental question in Reeb dynamics is the existence  of periodic orbits. A periodic orbit of $\lambda$ -- also called a closed orbit -- is a pair $P=(x,T), T>0,$ where $x:\R \to M$ is a periodic trajectory of the Reeb flow of $\lambda$ and $T$ is a period of $x$. We say that $P=(x,T)$ is simple if $T$ is the least positive period of $x$.


If $P=(x,T)$ is a periodic trajectory, then $x_T:= x(T\cdot):\R / \Z \to M$ is a critical point of the action functional $\mathcal{A}:C^\infty(\R / \Z,M) \to \R$, defined by $\gamma \mapsto \int_{\R / \Z} \gamma^* \lambda$. 
The critical value at $P$ coincides with the period $T$, see more in Appendix \ref{apendix_action_functional}. 


A closed orbit  $P=(x,T)$ is nondegenerate if the symplectic map
$D\psi_T: \xi|_{x(0)} \to \xi|_{x(0)}$ does not admit $1$ as an eigenvalue. Otherwise, $P$ is  degenerate. We say that $\lambda$ is nondegenerate if every periodic trajectory is nondegenerate. 

A  symplectic trivialization $\tau$ of $x^*\xi$ determines the generalized Conley-Zehnder index of a periodic orbit $P=(x,T)$,  denoted $\cz^\tau(P) \in \Z.$ 
This index is a measure of nearby orbits' rotation around $P$ along the period $T$. When available, we compute Conley-Zehnder indices in a fixed global trivialization of the contact structure. See Appendix \ref{appendix_CZrotation}. 

\subsection{Pseudo-holomorphic curves} If $\lambda$ is a contact form on a closed $3$-manifold $M$, then $(\R \times M, d(e^a\lambda))$ is a symplectic manifold, called the symplectization of $M$. Here, $a$ is the $\R$-coordinate.

Let $\xi=\ker \lambda$ be the contact structure and let $R$ be the Reeb vector field of $\lambda$. A complex structure $J:\xi \to \xi$ is called $d\lambda$-compatible if $d\lambda(\cdot, J \cdot)$ is a positive-definite inner product on $\xi$. 
Such a $J$ always exists and induces an $\R$-invariant almost complex structure $\widetilde J$ on $\R \times M$ satisfying
$$
\widetilde J|_\xi = J \ \ \mbox{ and } \ \ \widetilde J \cdot \partial_a = R.
$$

Let $(\Sigma,j)$ be a closed connected Riemann surface and let $\Gamma \subset \Sigma$
be a finite set of punctures. Denote $\dot \Sigma = \Sigma \setminus \Gamma$. A pseudo-holomorphic curve (or $\widetilde J$-holomorphic curve) in the symplectization is a map
$\tilde u=(a,u):\dot \Sigma  \to \R \times M$ 
satisfying
\begin{equation}\label{eq_phc}
\widetilde J(\tilde u) \circ d \tilde u = d\tilde u \circ j.
\end{equation}
This is a non-linear Cauchy-Riemann-type equation. 
In local complex coordinates $s+it \in \Sigma$,  the equation \eqref{eq_phc} assumes the form
$$
\pi u_s + J(u) \pi u_t =0, \ \ \ a_s = \lambda(u_t), \ \ \ a_t = -\lambda(u_s),
$$
where $\pi:TM \to \xi$ is the projection along the Reeb vector field.

The $d\lambda$-area of $\tilde u$ is 
the  integral $ \int_{\dot \Sigma} u^* d\lambda\geq 0,$
while the Hofer energy of $\tilde u$ is 
\begin{equation}\label{energiafinita} E(\tilde u) = \sup_{\psi
\in \Lambda} \int_{\dot \Sigma} \tilde u^* d(\psi(a) \lambda)\geq d\lambda\mbox{-area}(\tilde u).
\end{equation} Here, $\Lambda$ is the set of smooth functions $\psi:\R \to [0,1]$  satisfying $\psi'\geq 0$.   
We say that $\tilde u$ has finite energy if $0<E(\tilde u)<+\infty.$ In particular, a finite energy curve is not constant. 

The simplest example of a finite energy curve is the cylinder $$\tilde w(s,t) = (Ts, x_T(t)),\ \ \ \   \forall (s,t) \in \R \times \R / \Z,$$  over a periodic orbit $P=(x,T)$, where $x_T = x(T \cdot)$. The $d\lambda$-area of $\tilde w$ vanishes and its Hofer energy is the period $T>0$.



A puncture $z_0\in \Gamma$ of $\tilde u$ is called removable if $\tilde u$  smoothly extends to $\dot \Sigma \cup \{z_0\}$. Otherwise, $z_0$ is called non-removable. 
We assume that $z_0$ is non-removable and consider complex polar coordinates $(s,t) \in [0,+\infty)\times \R / \Z$
in a punctured neighborhood of $z_0\in \Sigma$.

The following fundamental theorem asserts that finite energy curves are asymptotic to periodic orbits at non-removable punctures.  

\begin{thm}[Hofer \cite{93}]\label{Ho93} If $z_0\in \Gamma$ is a non-removable puncture of $\tilde u=(a,u)$, then there exists $\epsilon \in \{-1,+1\}$, the sign of $z_0$, so that given any sequence $s_n \to +\infty$, there exists a
subsequence, still denoted by $s_n$, and a periodic orbit $P= (x,T),$
satisfying $u(s_n, \cdot) \to x(\epsilon T\cdot)$ in $C^\infty(\R / \Z,
M)$ as $s_n \to +\infty$.
\end{thm}

A non-removable puncture $z_0\in \Gamma$ is called positive or negative according to its sign given in Theorem \ref{Ho93}. The signs of the punctures induce the splitting $\Gamma = \Gamma^+ \cup \Gamma^-$. It can be shown that $a(z)\to  \pm\infty$  as $z\to z_0\in \Gamma^\pm$.  Since $\tilde u$ is not constant, Stokes' theorem implies that $\Gamma^+\neq \emptyset$.  

The orbit $P$ in Theorem \ref{Ho93} is called an asymptotic limit of $\tilde u$ at $z_0$. This asymptotic limit is unique when $P$ is nondegenerate \cite{props1}.  The curve $\tilde u$ may have multiple asymptotic limits at the same puncture. In this case, they are all degenerate, forming a continuum of periodic orbits with the same period.  See Siefring's example in \cite{Si3}.

The asymptotic operator $A_{P,J}$ associated with the closed orbit $P=(x,T)$ and the complex structure $J:\xi \to \xi$  acts on  sections of $x_T^* \xi$. Its eigenvalues and eigenfunctions determine the behavior of $\tilde u=(a,u)$ near a puncture. See Appendix \ref{appendix_asymptotic_operator} for more details.

\begin{figure}[ht]
\centering
\includegraphics[width=0.1\textwidth]{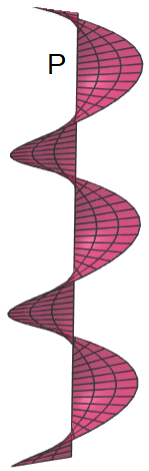}
\caption{Behavior of $u$ near a nondegenerate asymptotic limit $P$.}
\label{fig_winding}
\end{figure}

Assume that $\tilde u=(a,u)$ is asymptotic to a nondegenerate closed orbit $P_{z_0}$ at $z_0\in \Gamma^+$. Assume, moreover, that the $d\lambda$-area of $\tilde u$ is positive near $z_0$.
Then there exist an eigenvalue $\alpha< 0$ of $A_{P_{z_0},J}$ and an $\alpha$-eigenfunction $0\neq e(t)\in \xi|_{x_T(t)},\forall  t\in \R/ \Z,$  so that, for suitable coordinates $(s,t)\in [0,+\infty) \times \R / \Z$ near $z_0$,
$$
u(s,t) = {\rm exp}_{x_T(t)} \left\{e^{\alpha s}(e(t) + R(s,t))\right\}, \ \ \ \ \ \forall s \gg 0.
$$
The remainder term $R(s,t)$ and all of its derivatives decay exponentially fast to $0$ as $s\to +\infty$. In particular, $u$ does not intersect $P_{z_0}$ for every $s>0$ sufficiently large. If $z_0\in \Gamma^-$, then a similar asymptotic formula holds with a positive leading eigenvalue. See Figure \ref{fig_winding}.

In a suitable frame of $\xi$ along $P$, the winding number associated with the eigenfunction $e$, denoted by $\wind_\infty(\tilde u,z_0)\in \Z$,  satisfies
$$
\begin{aligned}
z_0 \mbox{ is positive } \Rightarrow \ \wind_\infty(\tilde u,z_0) \leq \frac{\cz(P_{z_0})}{2},\\
z_0 \mbox{ is negative } \Rightarrow\  \wind_\infty(\tilde u,z_0) \geq \frac{\cz(P_{z_0})}{2}.\\
\end{aligned}
$$

The Cauchy-Riemann nature of the pseudo-holomorphic curves implies that
$$
0 \leq \wind_\pi(\tilde u):= \sum_{z_0\in \Gamma^+}\wind_\infty(\tilde u,z_0) - \sum_{z_0\in \Gamma^-}\wind_\infty(\tilde u,z_0) - \chi(\Sigma)  + \#\Gamma,
$$ where $\chi(\Sigma)$ is the Euler characteristic of $\Sigma$ and the winding numbers are computed in a trivialization of the contact structure along the asymptotic limits  that  continuously extends to $u^*\xi$. 

The integer $\wind_\pi(\tilde u)$ counts (with multiplicity) the zeros of $\pi \circ du$, i.e., the number of points in $\dot \Sigma$ where $u$ is not transverse to the Reeb vector field. There are only finitely many such zeros, and each contributes positively.   

For instance, if $\tilde u=(a,u)$ is a pseudo-holomorphic plane (i.e. $\Sigma = \CP^1$ and $\Gamma = \{\infty\}$),  then
$$
0 \leq \wind_\pi(\tilde u) = \wind_\infty(\tilde u,\infty) -1 \leq \frac{\cz(P_\infty)}{2} -1,
$$
where $P_\infty$ is the asymptotic limit of $\tilde u$ at $\infty$. If $\cz(P_\infty)\in \{2,3\}$, then $\wind_\infty(\tilde u,\infty)=1$ and $\wind_\pi(\tilde u)=0$, and hence $u$ is everywhere transverse to $R$.

Now if $\tilde u=(a,u)$ is a pseudo-holomorphic cylinder ($\Sigma = \CP^1$ and $\Gamma = \{0,\infty\})$, with a positive puncture at $\infty$ and a negative puncture at $0$, then 
$$
0 \leq \wind_\pi(\tilde u) = \wind_\infty(\tilde u,\infty)-\wind_\infty(\tilde u,0) \leq  \frac{\cz(P_\infty)}{2}- \frac{\cz(P_0)}{2},
$$
where $P_\infty$  and $P_0$ are the asymptotic limits of $\tilde u$ at $\infty$ and $0$, respectively. If $\cz(P_\infty)=3$ and $\cz(P_0)=2$,  then $\wind_\infty(\tilde u,\infty)=\wind_\infty(\tilde u,0)=1$ and hence $\wind_\pi(\tilde u)=0$, which implies that $u$ is everywhere transverse to $R$.

The Fredholm index of $\tilde u$ is defined as
$$
{\rm Fred}(\tilde u)= \sum_{z\in \Gamma^+}\cz(P_z) - \sum_{z\in \Gamma^-} \cz(P_z) -\chi(\Sigma) + \# \Gamma,
$$
where $P_z$ is the asymptotic limit of $\tilde u$ at $z\in \Gamma$. The Conley-Zehnder indices are computed in a trivialization along the asymptotic limits that extends to $u^*\xi$. Under some favorable regularity conditions, ${\rm Fred}(\tilde u)$ gives the dimension about $\tilde u$ of the local family of pseudo-holomorphic curves with the same asymptotic data. 

If $\tilde u=(a,u)$ is a pseudo-holomorphic plane asymptotic to an index-$3$ orbit $P$, then ${\rm Fred}(\tilde u) = 2.$ The local two-dimensional family of planes asymptotic to $P$ contains $\R$-translations of $\tilde u$ and also displacements of $u$ in the direction of the flow, foliating a neighborhood of $\tilde u(\C)$.  If $\cz(P)=2$, and hence ${\rm Fred}(\tilde u) = 1$, then the local family consists of $\R$-translations of $\tilde u$ and, in this case, $\tilde u$ is called a rigid plane.

Finally, if $\tilde u$ is a pseudo-holomorphic cylinder asymptotic to an index-$3$ orbit at the positive puncture and to an index-$2$ orbit at the negative puncture, then ${\rm Fred}(\tilde u) = 1$. Hence the local family of cylinders with the same asymptotic data consists of $\R$-translations of $\tilde u$, and $\tilde u$ is called a rigid cylinder.

The pseudo-holomorphic planes and cylinders considered above are the building blocks of open book decompositions and, more generally, weakly convex  transverse foliations, discussed in the next sections. These foliations are projections of finite energy foliations in the symplectization.

\subsection{Global surfaces of section} Let $\psi_t:M \to M, t\in \R,$ be the Reeb flow of a contact form $\lambda$ on a smooth  closed  $3$-manifold $M$. We say that a compact embedded surface $\Sigma \hookrightarrow M$ is a {\bf global surface of section}  if
\begin{itemize}
\item $\partial \Sigma\neq \emptyset$ is a finite set of closed Reeb orbits;
\item $\Sigma \setminus \partial \Sigma$ is transverse to the flow;
\item every trajectory through a point in $M \setminus \partial \Sigma$ hits $\Sigma \setminus \partial \Sigma$  in the future and in the past.
\end{itemize}

The first return map $\Psi:\Sigma \setminus \partial \Sigma \to \Sigma \setminus \partial \Sigma$ preserves $d\lambda$ and encodes the dynamics of the Reeb flow. 

As an example, consider complex coordinates $(z_1=x_1+iy_1,z_2=x_2+iy_2) \in \C^2 \simeq \R^4.$ The Reeb flow of the standard contact form $\lambda_0$ on  $S^3$ is determined by
$$
\dot z_1 = 2i z_1, \ \ \ \ \ \dot z_2 = 2i z_2.
$$
The trajectories traverse the Hopf fibers  $z_1=ae^{2it},z_2=be^{2it},\forall t,$ where $(a,b)\in S^3.$  They are $\pi$-periodic, degenerate,  and have Conley-Zehnder index  $3$.

The standard contact structure $\xi_0=\ker\lambda_0$ is spanned by  the vector fields  on $S^3$
$$
\begin{aligned}
X_1 & = y_2\partial_{x_1} + x_2\partial_{y_1} -y_1 \partial_{x_2} - x_1 \partial_{y_2}, \\ X_2 & = -x_2 \partial_{x_1}+ y_2 \partial_{y_1} + x_1 \partial_{x_2} - y_1 \partial_{y_2},
\end{aligned}
$$
where $(x_1,y_1,x_2,y_2)\in S^3$. The vector fields $X_1$ and $X_2$ induce a global $d\lambda_0$-symplectic trivialization of $\xi_0$
$$
(aX_1+bX_2,d\lambda_0|_{\xi_0}) \mapsto ((a,b), dx \wedge dy).
$$

Let $J:\xi_0 \to \xi_0$ be the complex structure satisfying  $$J \cdot X_1=X_2.$$ 
For every $c\in S^1 \subset \C$, the map $$\tilde u_c(z) =(a_c(z),u_c(z))= \left(\frac{1}{4}\ln(|z|^2+1), \frac{(z,c)}{\sqrt{|z|^2+1}} \right), \ \ \ \forall z\in \C,$$
is a finite energy $\widetilde J$-holomorphic plane asymptotic to the Hopf fiber $P_1=S^1 \times \{0\}\subset \C^2$.  The closure of each $u_c(\C)$ is a global surface of section for the Reeb flow of $\lambda_0$, and  the family $u_c(\C),c\in S^1,$ spans an {\bf open book decomposition} whose binding is $P_1$.



If a hypersurface $i:S\hookrightarrow \R^4$ is  star-shaped  with respect to $0\in \R^4$, then the Reeb flow of $i^*\Lambda_0$ on $S$ conjugates with the Reeb flow of $\lambda= f\lambda_0$ on $S^3$, where $\Lambda_0$ is the Liouville form \eqref{eq_liou}, $\lambda_0$ is the standard contact form on $S^3$, and $f:S^3 \to (0,+\infty)$ is the unique function satisfying $\sqrt{f(p)}p\in S,\forall p\in S^3.$ 



\begin{defn}[Hofer-Wysocki-Zehnder]A contact form $\lambda=f\lambda_0$ on $S^3$ is called {\bf dynamically convex} if the Conley-Zehnder index of every closed trajectory  is  $\geq 3$.  \end{defn}

This terminology follows from the non-trivial fact that $i^*\Lambda_0$ is dynamically convex if the hypersurface $i:S\hookrightarrow \R^4$ is strictly convex, see \cite{convex}.

\begin{thm}[Hofer, Wysocki, Zehnder \cite{convex}]\label{th_convex}If $\lambda=f\lambda_0$ is a dynamically convex contact form on $S^3$, then its Reeb flow admits a disk-like global surface of section bounded by an index-$3$ closed Reeb orbit $P$.  \end{thm}


A remarkable corollary of Theorem \ref{th_convex} is that the Reeb flow of a dynamically convex contact form on $S^3$ has either $2$ or infinitely many periodic orbits. This follows from Brouwer's translation theorem and a result of J. Franks \cite{franks} on area-preserving homeomorphisms of the open annulus. The orbit $P$ in Theorem \ref{th_convex} is called a  {\bf Hopf fiber} since it is transversely isotopic to a Hopf fiber in $(S^3,\xi_0)$.



The disk-like global surface of section in Theorem \ref{th_convex} is the projection of a pseudo-holomorphic curve $$\tilde u_1=(a_1,u_1):\C \to \R \times S^3.$$ 
To find such a curve, one may start with a particular contact form $\lambda_0$ and complex structure $J_0$   for which a pseudo-holomorphic curve
$$
\tilde u_0=(a_0,u_0):\C \to \R \times S^3,
$$
is known to exist. Then,  interpolating the initial data  to the desired data, one tries to deform the existing curve $\tilde u_0$ through a family of pseudo-holomorphic curves
$$
\tilde u_\tau=(a_\tau,u_\tau):\C \to \R \times S^3, \ \ 0\leq \tau \leq 1,
$$
aiming at finding the desired curve $\tilde u_1$ at $\tau=1$.

Two crucial ingredients make this deformation mechanism work. Firstly, one needs an implicit function theorem that allows one to locally deform the curves  towards the target. The Fredholm theory in \cite{props3} shows that this deformation exists under appropriate circumstances, see also \cite{Wen1}. Secondly, one also needs a compactness theory describing the ends of these local families to ensure that the family eventually reaches the target.  The SFT-compactness theorem \cite{Abbas, BEHWZ} fulfills this part with precision.  

In Theorem \ref{th_convex}, the authors consider a symplectic cobordism from a suitable irrational ellipsoid to a hypersurface $S_\lambda$ corresponding to $\lambda$. Pseudo-holomorphic planes associated with the ellipsoid are then deformed towards $S_\lambda$ and, after an appropriate rescaling, they converge in the limit to a pseudo-holomorphic plane $\tilde u_0$, which is asymptotic to an index-$3$ orbit $P$ of the Reeb flow of $\lambda$. This orbit is expected to be the boundary of a disk-like global surface of section. To prove it,  $\tilde u_0$ is then deformed to establish an open book decomposition with binding orbit $P$.

\begin{figure}[ht!]
\centering
\includegraphics[width=0.5\textwidth]{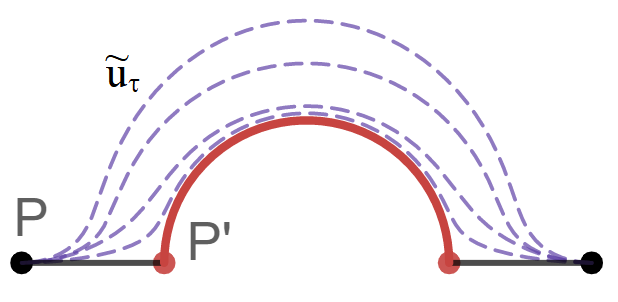}
\caption{The family of planes  asymptotic to $P$ is obstructed by the barrier formed by a rigid plane asymptotic to $P'$ and a cylinder from $P$ to $P'$.}
\label{fig_quebra}
\end{figure}

In the argument above, dynamical convexity is the main source of compactness. It is crucial that $\lambda$ does not have periodic orbits with low indices, especially those with Conley-Zehnder index $2$. The reason is that some non-deformable curves asymptotic to an index-$2$ orbit $P'$ might obstruct the continuation of the family of planes asymptotic to $P$, interrupting the construction of the desired open book decomposition, see Figure \ref{fig_quebra}. 

To avoid this type of obstruction, one usually imposes some geometric and topological assumptions on the Reeb flow such as linking conditions on the index-$2$ orbits. The following theorem provides conditions for a closed orbit $P$ to bound a disk-like global surface of section.

\begin{theo}[Hryniewicz, Salom\~ao \cite{HS2}]\label{thm_L_connected}
Let $\lambda=f\lambda_0$ be a contact form on $S^3$ and let $P \subset S^3$ be a closed Reeb orbit. Assume that $P$ is a Hopf fiber in $(S^3,\xi_0)$ and that the following hold:
\begin{itemize}
\item[(a)] $\cz(P)\geq 3$.
\item[(b)] every closed Reeb orbit $P'$ in the complement of $P$ satisfying either $\cz(P')=2$, or $\cz(P')=1$ and $P'$ is degenerate, is linked with $P$.
\end{itemize}
Then $P$ bounds a disk-like global surface of section.
\end{theo}

\begin{rem} The condition $\cz(P)\geq 3$ is not really necessary for $P$ to bound a disk-like global surface of section. Indeed, consider the Hamiltonian $$H = |z_1|^2 + |z_2|^4, \ \ \ \ \ \forall (z_1,z_2)\in \C^2.$$ The Reeb trajectories in $H^{-1}(1)$  are reparametrizations of $$z_1=ae^{it}, \ \ \ \ \ z_2=be^{2|z_2|^2i(t+c)}, \ \forall t,$$ where $a,b,c>0$ satisfy $a^2+b^4=1$.

The simple closed orbit $P_1= \{|z_1|=1,z_2=0\}$ has Conley-Zehnder index $1$. In spite of that, it bounds the disk-like global surface of section $\Sigma_1 = \{(z_1,z_2)\in H^{-1}(1): \im(z_2)=0 \mbox{ and } \re(z_2)  \in [0,1]\}$. In this case, however, the first return map does not extend to the closure of $\Sigma_1$.

It should be noticed that under a non-degeneracy assumption the condition $\cz(P) \geq 3$ is in fact necessary, see \cite{HS1}.

\end{rem}




Theorem \ref{thm_L_connected} generalizes the results in  \cite{char1,char2, hryn2,HS1}. Global surfaces of section with positive genus are studied in \cite{HSS}. Examples of Reeb flows on $(S^3,\xi_0)$ without disk-like global surfaces of sections are found by O. van Koert in \cite{vK,vk2}. Rational disk-like global surfaces of section on lens spaces are constructed in \cite{HLS,HS2,Sch}.


\subsection{Annulus-like global surfaces of section} G. Birkhoff \cite{Bi} showed the existence of annulus-like global surfaces of section for geodesic flows on positively curved $2$-spheres. H. Poincar\'e \cite{Po} showed that annulus-like global surfaces of section exist in the circular planar restricted three-body problem (CPR3BP) for energies below the first Lagrange value and small mass ratios. Both dynamical systems are Reeb flows on $$\R P^3 \equiv T^1S^2\equiv S^3 / \{\pm 1\},$$ equipped with the contact structure induced by $\xi_0=\ker \lambda_0$, see \cite{AFvKP}. 

The retrograde and direct orbits in the CPR3BP are $2:1$ projections of Hopf fibers and form a transverse link $L=L_1\cup L_2$ in  $(\R P^3,\xi_0)$, called a Hopf link. A link in $(\R P^3,\xi_0)$ that is transversely isotopic to $L$ is a Hopf link. A Hopf link binds an open book decomposition with annulus-like pages. If $\gamma$ is a simple closed geodesic on a Riemannian $2$-sphere, then the lift of $\dot \gamma$ and $-\dot \gamma$ to the unit sphere bundle $(\RP^3,\xi_0)$ forms a Hopf link $L$. The Birkhoff annulus $\Sigma \subset \R P^3$, formed by the unit vectors along $\gamma$  pointing to one of the hemispheres, is a page of an open book decomposition with binding $L$. 


The following theorem is a version of Theorem \ref{thm_L_connected} for Reeb flows on $(\R P^3,\xi_0)$. It provides conditions for a pair of closed Reeb orbits forming a Hopf link to bound an annulus-like global surface of section.

\begin{theo}[Hryniewicz, Salom\~ao, Wysocki \cite{HSW}]\label{main2}
Let $\lambda=f\lambda_0$ be a contact form on $(\R P^3,\xi_0)$ and let $L=L_1\cup L_2$ be a Hopf link formed by closed Reeb orbits. Assume
that
\begin{itemize}
\item[(a)] $\cz(L_i)>0, i=1,2.$ 
\item[(b)] all periodic orbits in $\R P^3\setminus L$ satisfying $\cz\in \{-3,\ldots,5\}$ have non-zero intersection number with a Birkhoff annulus bounded by $L$.
\end{itemize}
Then  $L$ bounds an annulus-like global surface of section for the Reeb flow.
\end{theo}

 In \cite{HSW} one finds a generalization of Theorem \ref{main2} for planar contact $3$-manifolds. The source of pseudo-holomorphic curves in this case is a particular choice of contact form and complex structure admitting a planar open book decomposition whose pages are projections of pseudo-holomorphic curves, see Chris Wendl's construction in \cite{Wen2}. In a suitable symplectic cobordism, one deforms these  pseudo-holomorphic curves to obtain an open book decomposition adapted to the Reeb flow of $\lambda$. 
 
 The linking assumptions in Theorem \ref{main2}-(b) are the main source of compactness, ruling out potential  rigid curves that may obstruct the construction of the open book decomposition.

Our next task is to consider transverse foliations which are not necessarily open books. Such foliations naturally rise due to the lack of compactness of a family of pseudo-holomorphic curves.

\section{Transverse foliations} Transverse foliations  generalize the concept of open book decomposition adapted to the flow. They may be used as tools to address relevant dynamical questions such as the existence and multiplicity of closed orbits, homoclinics, heteroclinics, topological entropy, etc.

\begin{defi}\label{def_transverse}
Let $(M,\xi)$ be a smooth closed connected contact $3$-manifold and let $\lambda$ be a defining contact form, i.e. $\ker \lambda=\xi$. A transverse foliation of $M$ adapted to the Reeb flow of $\lambda$ is a singular foliation $\mathcal{F}$ of $M$ so that:
\begin{itemize}

\item[(i)] the singular set  $\mathcal{S}\subset \mathcal{F}$ is formed by finitely many simple closed Reeb orbits $P_1,\ldots ,P_l$, called  binding orbits. The link $L=\bigcup_i P_i$ is the binding of $\mathcal{F}$.

\item[(ii)] a leaf of $\mathcal{F} \setminus \mathcal{S}$  is a properly embedded surface $\dot \Sigma \hookrightarrow M \setminus L$. Its closure is a compact embedding  $\Sigma\hookrightarrow M$ with non-empty boundary, and whose boundary components are binding orbits.  

\item[(iii)] every $\dot \Sigma$ is transverse to the flow. The orientation of $M$ ($\lambda \wedge d\lambda>0$) and the Reeb vector field orient $\dot \Sigma$, which orients $\partial \Sigma$ as a boundary. A component of $\partial \Sigma$ is called positive if its orientation agrees with the one induced by the flow.  Otherwise, it is called negative. Each end of $\dot \Sigma$ is called a puncture and the corresponding periodic orbit in $\partial \Sigma$ is called the asymptotic limit of $\dot \Sigma$ at the puncture.


\end{itemize}
\end{defi}

The following fundamental result states that a generic Reeb flow on the tight $3$-sphere admits a transverse foliation.

\begin{theo}[Hofer-Wysocki-Zehnder \cite{fols}]\label{sfef}
Let $\lambda=f\lambda_0$ be a nondegenerate contact form on $(S^3,\xi_0)$. Then the Reeb flow of $\lambda$ admits a transverse foliation.  The binding orbits have Conley-Zehnder indices $1$, $2$ or $3$ and the regular leaves have genus zero. Every regular leaf  has precisely one positive puncture and an arbitrary number of negative punctures. One of the following cases holds:
\begin{itemize}
\item the asymptotic limit of $\dot \Sigma$ at its positive puncture has index $3$ and the asymptotic limit of $\dot \Sigma$ at every negative puncture has index $1$ or $2$. There exists at most one negative puncture whose asymptotic limit has index $2$.
\item the asymptotic limit of $\dot \Sigma$ at its positive puncture has index $2$ and the asymptotic limit of $\dot \Sigma$ at every negative puncture has index $1$.
\end{itemize}
\end{theo}

The transverse foliation $\F$ in Theorem \ref{sfef} is the projection to $S^3$ of a finite energy foliation, that is an  $\R$-invariant foliation of $\R \times S^3$ whose leaves are the image of finite energy punctured spheres $\tilde u=(a,u):S^2 \setminus \Gamma \to \R \times S^3$. The source of these pseudo-holomorphic curves is a family of pseudo-holomorphic spheres on $\CP^2$. 

The simplest example of transverse foliation is an open book decomposition in $S^3$ with precisely one binding orbit $P_3$ having Conley-Zehnder index 3. In this case, $S^3 \setminus P_3$ is foliated by an $S^1$- family of planes asymptotic to $P_3$ and each regular leaf is a disk-like global surface of section for the Reeb flow.

Now let us assume that $\F$ is a transverse foliation such that every binding orbit has index $2$ or $3$. Denote by $P_{2,i},i=1,\ldots,l,$ the index-$2$ binding orbits and by $P_{3,j},j=1,\ldots,l+1,$ the index-$3$ binding orbits. In this case $\F$ is called a weakly convex foliation. They are organized as follows:
\begin{itemize}
    \item  each $P_{2,i}$ is the asymptotic limit of two rigid planes $U_{i,1}$ and $U_{i,2}$, so that $\S_i = U_{i,1} \cup P_{2,i} \cup U_{i,2}$ is a $C^1$-embedded $2$-sphere.
    
    \item  each component $\U_j \subset S^3 \setminus \bigcup_{i=1}^l \S_i,j=1,\ldots,l+1,$ contains a binding orbit $P_{3,j}$. Given $\S_i\subset \partial \U_j$ there exists a  rigid cylinder $V^j_i\subset \U_j$ asymptotic to $P_{3,j}$ at its positive puncture and to $P_{2,i}$ at its negative puncture.
    
    \item  if $k_j\in \N^*$ is the number of boundary components of $\U_j$, then there exist precisely $k_j$ one-parameter families of planes asymptotic to $P_{3,j}$. At its ends, such a family breaks onto a rigid cylinder and a rigid plane. 
\end{itemize}


\begin{defi}A contact form $\lambda=f\lambda_0$ on $(S^3,\xi_0)$  is called {\bf weakly convex} if all of its closed trajectories have Conley-Zehnder index $\geq 2$.
\end{defi}
 
 Theorem \ref{sfef} says that if $\lambda=f\lambda_0$ is a nondegenerate weakly convex contact form on $(S^3,\xi_0)$, then its Reeb flow admits a weakly convex foliation. One particular weakly convex foliation is the 3-2-3 foliation.





\begin{figure}[ht]
\centering
\includegraphics[width=1.0\textwidth]{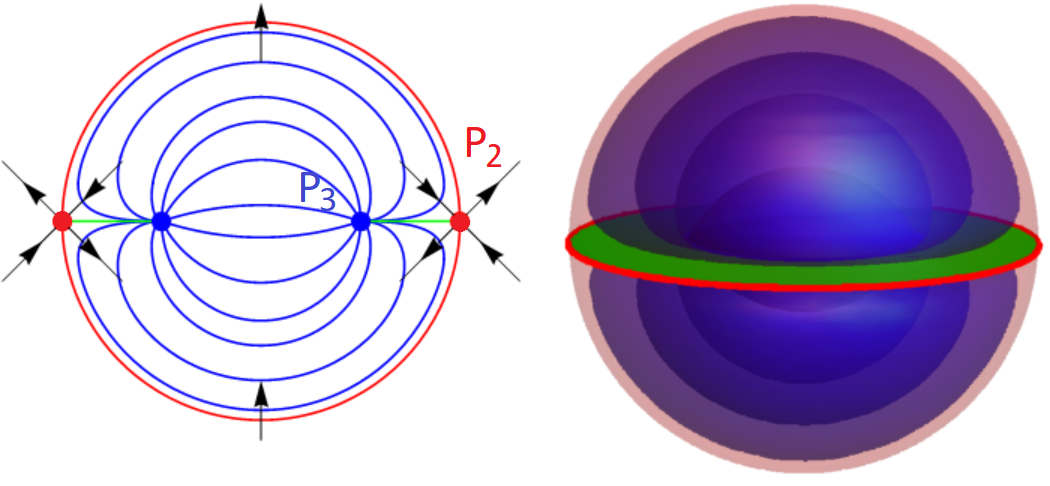}
\caption{The 3-2-3 foliation inside one of the $3$-balls.}
\label{fig_323}
\end{figure}

\begin{defi}
A {\bf 3-2-3 foliation} adapted to the Reeb flow of a contact form $\lambda=f\lambda_0$ on $(S^3,\xi_0)$ is a weakly convex foliation with precisely three binding orbits $P_3$, $P_2$ and $P_3'$. In particular, they are unknotted, mutually unlinked and their respective Conley-Zehnder indices are $3,2$, and $3$. The regular leaves are the following:
\begin{itemize}

\item a pair of rigid planes $U_1$ and $U_2$ asymptotic to $P_2$. The $2$-sphere $U_1 \cup P_2 \cup U_2$ separates $S^3$ into two components diffeomorphic to open $3$-balls $\mathcal{B}$ and $\mathcal{B}'$, so that the orbits $P_3\subset \mathcal{B}$ and $P_3'\subset \mathcal{B}'$. 

\item a cylinder $V\subset \mathcal{B}$ asymptotic to $P_3$ at its positive puncture and to $P_2$ at its negative puncture; and a cylinder $V'\subset \mathcal{B}'$ asymptotic to $P_3'$ at its positive puncture and to $P_2$ at its negative puncture.

\item a one-parameter family of planes in $\mathcal{B}$, asymptotic to $P_3$ at their positive punctures; and a one-parameter family of planes in $\mathcal{B}'$, asymptotic to $P_3'$ at their positive punctures.
\end{itemize}
\end{defi}


The 3-2-3 foliations are motivated by the study of Hamiltonian dynamics near certain critical energy surfaces.

\subsection{Critical energy surfaces} 

Periodic orbits with Conley-Zehnder index $2$ can only exist in the absence of convexity.  As an example, consider the Hamiltonian
$$
H_0 = -\frac{\alpha}{2}\im (z_1^2) + \frac{\omega}{2}|z_2|^2, \ \ \ \ (z_1,z_2) \in \C^2,
$$
where $\alpha,\omega>0$. The Reeb trajectories are reparametrizations of
$$
z_1 = ae^{-\alpha t} + i be^{\alpha t}, \ \ \ z_2 = ce^{-\omega i(t+d)}, \ \ \ \forall t,
$$
for $a,b,c ,d\in \R$. They behave like a saddle in the $z_1$-plane and like a center in the $z_2$-plane.

\begin{figure}[ht!]
\centering
\includegraphics[width=0.7\textwidth]{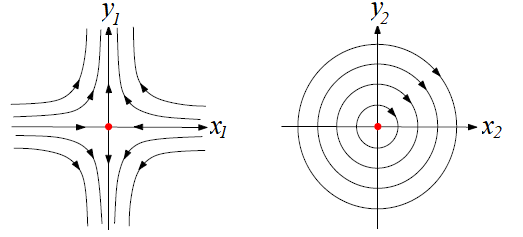}
\caption{The trajectories near a saddle-center. }
\label{fig_esfera0}
\end{figure}

The critical point  $0\in H^{-1}_0(0)$ is a rest point of the flow. For each $E>0$, the energy level $H^{-1}_0(E)$ contains an index-$2$ orbit
$$P_{2,E} = \left \{z_1=0, |z_2|^2 = \frac{2E}{\omega}\right \},$$ which bounds the pair of rigid planes
$$
\begin{aligned}
U_{1,E} =\left \{\im (z_1) = -\re (z_1)>0,|z_2|^2<\frac{2E}{\omega}\right \}\cap H^{-1}_0(E),\\  U_{2,E} =\left \{\im (z_1) = -\re (z_1)<0,|z_2|^2<\frac{2E}{\omega}\right \}\cap H^{-1}_0(E).
\end{aligned}
$$

Every $P_{2,E}$ is hyperbolic inside $H^{-1}_0(E)$ and its stable and unstable manifolds project to the real and imaginary axes of the $z_1$-plane.

The $2$-sphere $U_{1,E} \cup P_{2,E} \cup U_{2,E}$ separates  $H^{-1}_0(E)$ into two (unbounded) components, each one  of them contains a branch of the stable and a branch of the unstable manifold of $P_{2,E}$:
$$
\begin{aligned}
W^s(P_{2,E}) & = \left \{\im (z_1)=0, |z_2|^2=\frac{2E}{\omega} \right \},\\
W^u(P_{2,E}) & =\left \{\re (z_1)=0, |z_2|^2=\frac{2E}{\omega} \right \}.
\end{aligned}
$$

The dynamics of $H_0$ near the origin is a toy model for the dynamics near a saddle-center rest  point. 

\begin{defi} Let $p_c\in \R^4$ be a critical point of a smooth function $H: \R^4 \to \R$. We say that $p_c$ is a saddle-center equilibrium if the linearization of the Hamiltonian vector field at $p_c$  admits a pair of real eigenvalues and a pair of purely imaginary eigenvalues.
\end{defi}

To illustrate the existence of transverse foliations, assume that $H$ is a real-analytic Hamiltonian function admitting a saddle-center $p_c \in H^{-1}(0)$.

\begin{figure}[ht!]
\centering
\includegraphics[width=1.0\textwidth]{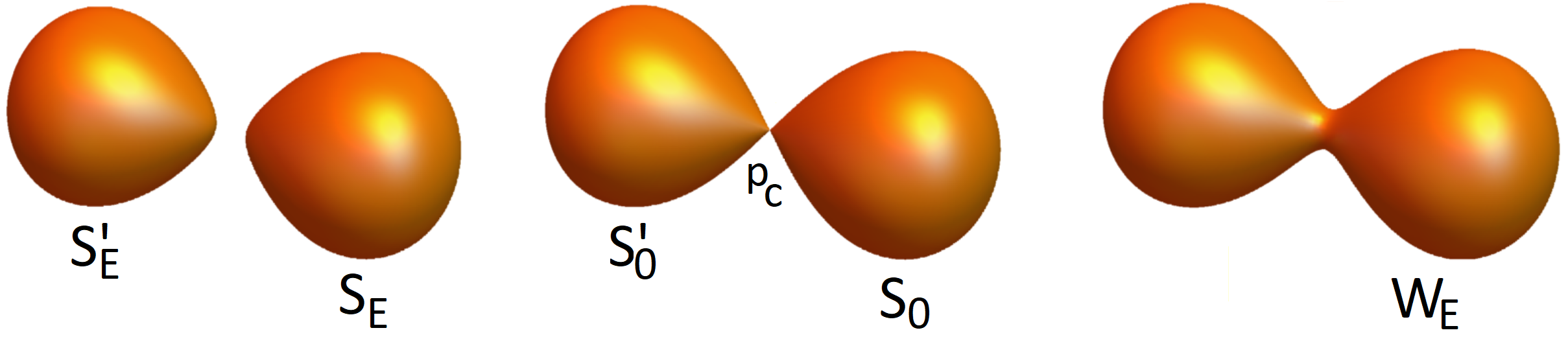}
\caption{The energy levels for $E<0, E=0$ and $E>0$.}
\label{fig_Niveis}
\end{figure}

Assume that 
\begin{itemize}
    \item  for $E<0$, $H^{-1}(E)$ contains two sphere-like components $S_E$ and $S_E'$, $C^0$-converging to $S_0$ and $S_0'$, respectively, as $E \to 0^-$. 
    \item   the subsets $S_0,S_0'\subset H^{-1}(0)$ intersect each other at $p_c$, and are called singular sphere-like hypersurfaces.
    \item for $E>0$,  $H^{-1}(E)$ has a sphere-like component $W_E$ corresponding to the connected sum $S_E \# S_E'$.
\end{itemize}
See Figure \ref{fig_Niveis}.

 The Hamiltonian $H$, in suitable coordinates  $(z_1,z_2)\in \C^2$ near $p_c\equiv 0,$ is assumed to have the form
$$
H=H_0 + R(\im(z_1^2),|z_2|^2),
$$
where $R(\cdot,\cdot)$ is a real-analytic function which vanishes up to first order at $(0,0)$. 

\begin{figure}[ht!]
\centering
\includegraphics[width=0.9\textwidth]{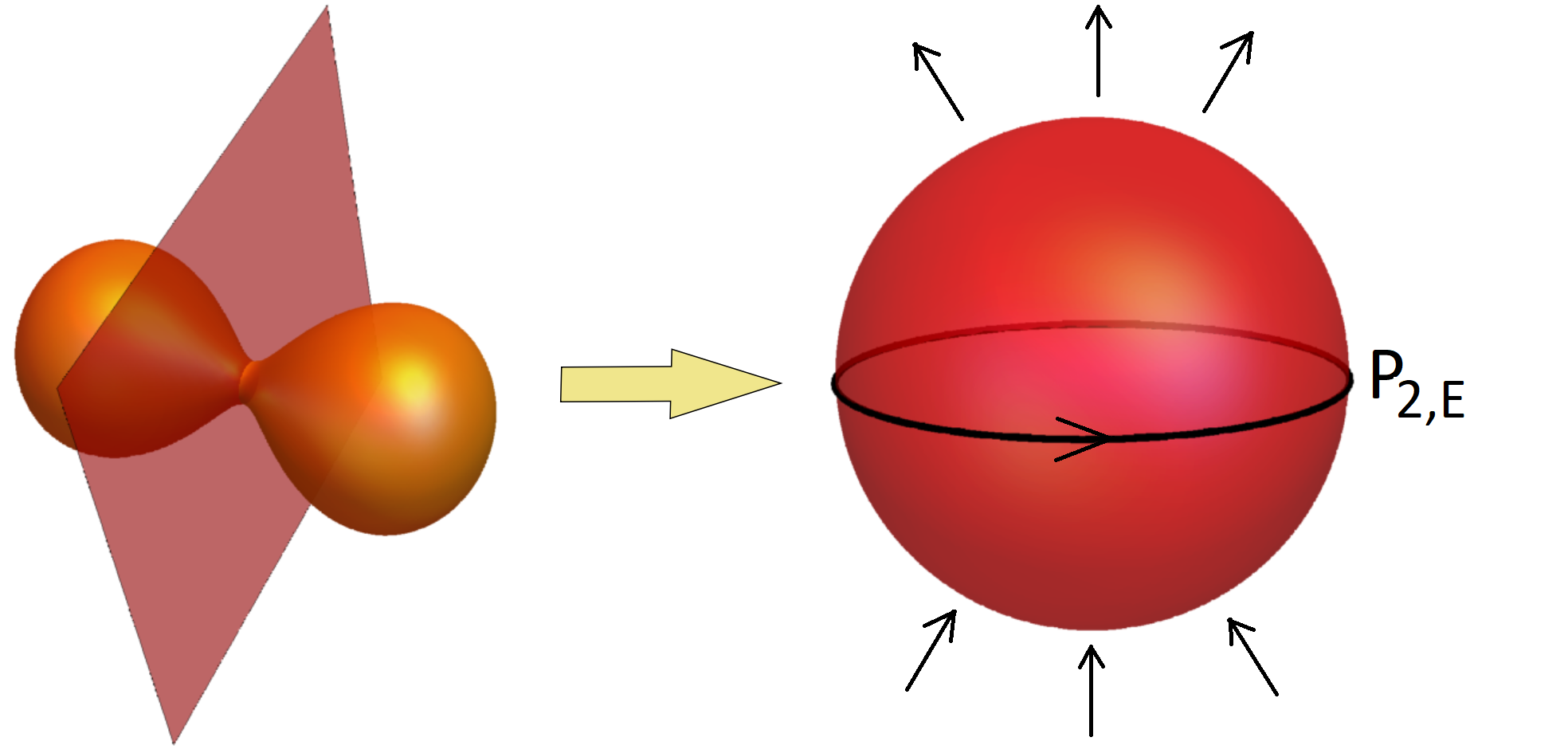}
\caption{The separating $2$-sphere 
in the neck region of $W_E$.}
\label{fig_corte}
\end{figure}

In the neck region of $W_E$, there exists an index-$2$ closed orbit $P_{2,E}$ -- which will be referred to as the Lyapunoff orbit -- and a pair of rigid planes $U_{1,E}$ and $U_{2,E}$ asymptotic to $P_{2,E}$, see Figure \ref{fig_corte}. The $2$-sphere $U_{1,E}\cup P_{2,E} \cup U_{2,E}$ separates $W_E$ into two open $3$-balls. The rigid planes are transverse to the flow and correspond to transit trajectories from one $3$-ball to the other.

\begin{defi}We say that $S_0\subset H^{-1}(0)$ containing a saddle-center equilibrium $p_c$ as its unique singularity is strictly convex if $S_0$ bounds a convex domain in $\R^4$ and all the sectional curvatures of $S_0\setminus \{p_c\}$ are positive.
\end{defi}

A similar definition holds for $S_0'$. We aim at finding a transverse foliation of $W_E$ so that $P_{2,E}$ is a binding orbit and the hemispheres $U_{1,E}$ and $U_{2,E}$ are leaves.  Actually, such a foliation exists provided that both $S_0$ and $S_0'$ are strictly convex.

\begin{theo}[de Paulo, Salom\~ao \cite{PS1,PS2}]\label{th_dPS} Assume that both the singular sphere-like hypersurfaces $S_0$ and $S_0'$ are strictly convex. Then, for every $E>0$ sufficiently small, the Hamiltonian flow on $W_E$ admits a $3-2-3$ foliation.  Moreover, $P_{2,E}$ is a binding orbit and the hemispheres $U_{1,E}$ and $U_{2,E}$ are leaves. In addition, there exist infinitely many periodic orbits and infinitely many homoclinics to $P_{2,E}$ in both sides of the separating $2$-sphere.
\end{theo}

Notice that there are no non-degeneracy assumptions in Theorem \ref{th_dPS}. Its proof involves some careful estimates of Conley-Zehnder indices of the periodic orbits in $W_E$. Using the quaternionic structure of $\R^4$, see \cite{GS1}, the transverse linearized flow is given by the solutions to the equation
$$
\left(\begin{array}{c}\dot  a_1 \\ \dot a_2 \\  \end{array} \right) =-J M\left(\begin{array}{c} a_1 \\ a_2 \end{array} \right).
$$
Here, $M$ is a symmetric matrix whose entries depend on the second order derivatives of $H$ along the Hamiltonian trajectory,
see Appendix \ref{ap_quaternios} for more details.
The strict convexity of $S_0$ and $S_0'$  implies that outside a small neighborhood $\U$ of $p_c$, the matrix $M$ is uniformly positive-definite, for every $E>0$ sufficiently small.  In particular, the solutions $a_1+ia_2 = re^{i\eta}\in \C$ rotate uniformly in the counter-clockwise direction, i.e., there exists $K_0>0$ so that $\dot \eta(t) > K_0, \forall t$. This gives a uniform estimate of the form $\cz(P)> \beta T$ for some $\beta>0$, as long as the periodic orbit $P=(x,T)$ does not intersect $\U$. However, a local estimate involving the linearized flow near a saddle-center shows that the contribution to the variation of $\eta$ tends to be proportional to the amount of time spent near $p_c$. Therefore, if $\U$ is taken sufficiently small, then the periodic orbits intersecting $\U$ spend a lot of time near $p_c$ and their Conley-Zehnder indices are arbitrarily high. Summarizing, the convexity of $S_0$ and $S_0'$ implies that the only trajectory of $W_E$ with index less than $3$ is the Lyapunoff orbit $P_{2,E}$, whose index is $2$, and the closed orbits passing close to $p_c$ tend to have arbitrarily high indices. In particular, the transit periodic orbits, i.e., those which cross the separating $2$-sphere, have arbitrarily high indices if $E>0$ is sufficiently small and, moreover, $P_{2,E}$ is not linked to any index-3 periodic orbit. 

The Hamiltonian flow on $W_E$ is equivalent to the Reeb flow of a contact form $\lambda_E$, which is weakly convex due to the estimates of the Conley-Zehnder indices mentioned above. Taking a sequence $\lambda_n\to \lambda_E$ in $C^\infty$ of nondegenerate contact forms on $W_E$, we show that the above estimates hold for $\lambda_n$, for every large $n$. Applying Theorem \ref{sfef}, we conclude that $\lambda_n$ admits a $3-2-3$ foliation so that the Lyapunoff orbit $P_{2,E}$ is a binding orbit. Now for particular choices of almost-complex structures, the bubbling-off analysis of pseudo-holomorphic curves allows one to push the foliation of $\lambda_n$ to a $3-2-3$ foliation adapted to the Reeb flow of $\lambda_E$ as in Theorem \ref{th_dPS}.  

The $3-2-3$ foliation organizes the flow in $W_E$. By thoroughly studying some transition maps defined by the flow between its leaves -- which are area-preserving surface maps --, one may obtain infinitely many homoclinics to $P_{2,E}$ and infinitely many periodic orbits in $W_E$.

As an example, assume that $H$ has the form kinetic plus potential energy
\begin{equation}\label{Hpot}H = \frac{y_1^2+y_2^2}{2} + V(x_1,x_2).
\end{equation}
If $x_c\in \R^2$ is a saddle-type critical point of $V$, then $p_c=(x_c,0)\in \R^4$ is a saddle-center equilibrium of $H$. 

Assume that $p_c \in H^{-1}(0)$ and that $H^{-1}(0)=S_0 \cup S_0'$, where $S_0,S_0'$ are singular sphere-like hypersurfaces intersecting at $p_c$. The projection of $H^{-1}(0)$ to the $(x_1,x_2)$-plane is the Hill's region $D_0 \cup D_0'$, where $D_0$ and $D_0'$ are topological disks with smooth boundaries except at $\{x_c\}=D_0 \cap D_0'$. The strict convexity of $S_0$ or $S_0'$ can be checked  according to the following theorem.

\begin{theo}[\cite{Sa1}] \label{th_potencial_convexo}The sphere-like component $S_0\subset H^{-1}(0)$ is strictly convex if and only if the inequality  $$-2V(V_{x_{1}x_{1}}V_{x_{2}x_{2}}-V_{x_{1}x_{2}}^{2})+V_{x_{1}x_{1}}V_{x_{2}}^{2}+
V_{x_{2}x_{2}}V_{x_{1}}^{2}-2V_{x_{1}}V_{x_{2}}V_{x_{1}x_{2}}>0,$$
holds on $D_0\setminus \{x_c\}$. A similar statement holds for $S_0'$.
\end{theo}

This criterion was used to check the convexity properties of several Hamiltonian systems, see \cite{PS2,Kim3,Sch}. See  \cite{GS1,GS2} for results on the dynamics on the singular sphere-like subsets of the critical energy surface.  Transverse foliations admitting index $1$-binding orbits are found in \cite{lemos,lemosth,Wen0}. Transverse foliations are stable under connected sums, as shown by Fish and Siefring in \cite{FS}. An interesting approach to transverse foliations using techniques from embedded contact homology are obtained by Colin, Dehornoy and Rechtman in \cite{CDR}.

\appendix
\section{The Action Functional }\label{apendix_action_functional}

Let $M$ be a smooth $3$-manifold equipped with a contact form $\lambda$. Denote by $\varphi_t,t\in \R,$ the Reeb flow of $\lambda$.

The action functional associated with $\lambda$ is defined on the space of smooth curves $\gamma: \R / \Z \to  M$ as 
$$
\mathcal{A}(\gamma) = \int_{\R / \Z} \gamma^*\lambda, \ \ \ \ \ \forall \gamma\in C^\infty(\R / \Z,M).
$$




Let $\eta \in \Gamma(\gamma^*TM)$ be a vector field along $\gamma$ and let $u:(-\epsilon,\epsilon)\times \R / \Z \to M$ be a variation of $\gamma$ so that $u(0,t)=\gamma(t)$ for all $t\in \R / \Z$ and $u_s(0,\cdot)=\eta$. We compute the first variation of $\mathcal{A}$ at $\gamma$ in the direction of $\eta$
\begin{equation}\label{1derivaA}
\begin{aligned}
d\mathcal{A}(\gamma) \cdot \eta
& = \int_{0 \times \R / \Z} \mathcal{L}_{\partial_s} (u^*\lambda)\\
& = \int_{0 \times \R / \Z} d i_{\partial_s} (u^*\lambda)+\int_{0 \times \R / \Z} i_{\partial_s} (u^*d\lambda)\\
& = \int_0^1 d\lambda|_{\gamma(t)}(\eta(t), \dot \gamma(t))dt.
\end{aligned}
\end{equation}

Since $d\lambda|_{\xi=\ker \lambda}$ is nondegenerate, this implies that $\gamma$ is a critical point of $\mathcal{A}$ if and only if $\dot \gamma \subset \ker d\lambda|_\gamma.$ In particular, if $\dot \gamma$ never vanishes, then $\gamma$ is a critical point of $\mathcal{A}$ if and only if $\gamma$ is a reparametrization of a closed Reeb orbit.

Now let $P=(x,T)$ be a closed Reeb orbit of $\lambda$. Then  $x_T = x(T\cdot):\R / \Z \to M$ is a critical point of $\mathcal{A}$. Let $\eta \in \Gamma(x_T^*\xi)$ and let $u:(-\epsilon,\epsilon) \times \R / \Z \to M$ be a variation of $x_T$ so that $u_s(0,\cdot)=\eta$. We can assume that $u_s(s,\cdot)\subset \xi, \forall s$. Using a finite  covering of a tubular neighborhood of $x_T(\R)$, we may assume that $x_T$ is an embedding. Let $\Xi$ be a vector field extending $\eta$ in a neighborhood of $x(\R)$. 

 The second variation of $\mathcal{A}$ at $x_T$ in the direction of $\eta$ is
\begin{equation}\label{second_derivative}
\begin{aligned}
d^2\mathcal{A}(x_T) \cdot (\eta,\eta) 
& = \int_{0 \times \R / \Z} \mathcal{L}_{\partial_s}(\mathcal{L}_{\partial_s} (u^*\lambda))\\
& = \int_{0 \times \R / \Z} i_{\partial_s}d (i_{\partial_s} (u^*d\lambda)))\\
& = \int_{ \R / \Z} x_T^* (i_\Xi d(i_\Xi d\lambda))\\
& = \int_0^1 d(i_\eta d\lambda)|_{x_T(t)}(\eta(t), \dot x_T(t))dt\\
& = \int_0^1 d\lambda|_{x_T(t)} ( -(\mathcal{L}_{\dot x_T} \eta)(t),\eta(t)) \, dt.\\
\end{aligned}
\end{equation}
We have used the usual relation $d\alpha(X,Y)=X \cdot \alpha(Y) - Y \cdot \alpha(X) - \alpha ([X,Y])$, where $X,Y$ are vector fields and $\alpha$ is a $1$-form.

Observe that the last expression obtained in \eqref{second_derivative} depends only on the linearized flow along $P$
$$
\begin{aligned}
(\mathcal{L}_{\dot x_T} \eta)(t) & : = \frac{d}{ds}\big|_{s=0} \left\{ D\varphi^{-1}_{sT}(x_T(t+s)) \cdot \eta(t+s) \right\}.
\end{aligned}
$$



A similar analysis of the action functional can be found in \cite{lemos}.

\section{The asymptotic operator}\label{appendix_asymptotic_operator} Let $P=(x,T)$  be a closed orbit of the Reeb flow of $\lambda$ and let $J:\xi \to \xi$ be a $d\lambda$-compatible complex structure. The asymptotic operator associated with $P$ and $J$ is the linear operator 
$$
\begin{aligned}
A_{P,J}(\eta) := -J|_{x_T} \cdot  \mathcal{L}_{\dot x_T}\eta\in L^2(x_T^* \xi), \ \ \ \ \ \ \ \forall \eta\in W^{1,2}(x_T^*\xi), \\
\end{aligned}
$$  
where $\mathcal{L}_{\dot x_T} \eta$ is defined in the previous section.

It is an exercise to check that  $A_{P,J}$ admits $0$ as an eigenvalue if and only if $P$ is degenerate.

Since $J$ preserves $d\lambda$ we see from \eqref{second_derivative} that
$$
d^2\mathcal{A}(x_T) \cdot (\eta,\eta) = \int_0^1 d\lambda|_{x_T(t)} (A_{P,J} \cdot \eta (t), J_t\cdot \eta(t))\, dt,
$$ where $J_t = J|_{x_T(t)}$. 

If $\nu\in\R$ is an eigenvalue of $A_{P,J}$ and $\eta_\nu$ is a non-trivial $\nu$-eigenfunction, then
$$
\begin{aligned}
d^2\mathcal{A}(x_T) \cdot (\eta_\nu,\eta_\nu) & =  \int_0^1 d\lambda_{x_T(t)}(A_{P,J} \cdot \eta_\nu (t), J_t\cdot \eta_\nu(t))\, dt,\\
& = \nu \int_0^1 g_{J_t}(\eta_\nu(t),\eta_{\nu}(t))\, dt,
\end{aligned}
$$ where $g_J(\cdot, \cdot) = d\lambda(\cdot, J \cdot)$ is the positive-definite inner product on $\xi$ induced by $d\lambda$ and $J$. Hence $d^2\mathcal{A}(x_T) \cdot (\eta_\nu,\eta_\nu)$ vanishes if and only if $\nu=0$ and, otherwise, its sign coincides with the sign of $\nu$. 

In order to better describe the spectrum of the operator $A_{P,J}$ we need some normalization.
Choose a unitary trivialization 
 $\Psi: x_T^* \xi \to \R / \Z \times \R^2$. This means that
$$
\Psi^*(dx \wedge dy) = d\lambda|_\xi \ \ \ \ \mbox{ and } \ \ \ \ \ \Psi_* J = J_0:=\left(\begin{array}{cc}0 & -1 \\ 1 & 0 \end{array}\right).
$$ 
Here $(x,y)$ are coordinates in $\R^2$. 


We claim that the operator 
$$
\mathcal{L}_S:= \Psi \circ A_{P,J} \circ \Psi^{-1}
$$ has the form
\begin{equation}\label{operator_unitary}
\mathcal{L}_S= -J_0\frac{d}{dt} - S,
\end{equation} 
 for some smooth loop $t \mapsto S(t)$ of $2\times 2$ symmetric matrices. 
 
 Indeed, recall that the Reeb flow of $\lambda$ preserves $d\lambda|_\xi$ and thus we find a smooth path of $2\times 2$ symplectic matrices $\Phi(t)$, satisfying  $\Phi(0)=I$ and
\begin{equation}\label{periodicity}
\Phi(1+t) = \Phi(t) \Phi(1), \ \ \ \ \ \  \forall t\in \R / \Z.
\end{equation}
It represents the linearized flow of $\frac{1}{T}\lambda$ restricted to $(\xi,d\lambda|_\xi)$ along $x_T$ in coordinates $\Psi$. In particular, a solution $\bar \zeta(t)\in \R^2 \simeq \xi|_{x_T(t)}$ to the linearized flow  with initial condition $\bar \zeta(t+s)=\zeta(t+s) \in \R^2 \simeq \xi_{x_T(t+s)}$ satisfies
$$
\bar \zeta(t)  = \Phi(t)\Phi(t+s)^{-1}\zeta(t+s), \ \ \ \ \ \  \forall s,t.
$$ Note that $\bar \zeta(t)$ depends on $s$.

Denoting $$\zeta(t) = \Psi \circ \eta(t)\in \R^2  \ \forall t,$$ let us compute $\mathcal{L}_{\dot x_T} \eta$ in coordinates $\Psi$ 
\begin{equation}\label{Liederivative}
\begin{aligned}
\Psi \circ \mathcal{L}_{\dot x_T}\eta (t) & = \lim_{s\to 0} \frac{\bar \zeta(t) - \zeta(t)}{s}\\
& = \lim_{s\to 0} \frac{\Phi(t)\Phi(t+s)^{-1}\zeta(t+s) - \zeta(t)}{s}\\
& = \lim_{s\to 0} \left\{\Phi(t)\Phi(t+s)^{-1}\frac{\zeta(t+s)-\zeta(t)}{s}+\Phi(t)\frac{\Phi(t+s)^{-1}-\Phi(t)^{-1}}{s}\zeta(t)\right\}\\
& = \dot \zeta(t) +\Phi(t)\dot{[\Phi(t)^{-1}]}\zeta(t)\\
& = \dot \zeta(t) - \dot \Phi(t) \Phi(t)^{-1}\zeta(t)
\end{aligned}
\end{equation}

Now let 
\begin{equation}\label{definitionS}
S(t):= -J_0 \dot \Phi(t) \Phi(t)^{-1}.
\end{equation} Since $\Phi(t)$ is symplectic  we have
$$
\Phi(t)^T J_0 \Phi(t) = J_0 \Rightarrow \dot \Phi(t)^T J_0 \Phi(t) + \Phi(t)^T J_0 \dot \Phi(t) = 0, \ \ \ \forall t.
$$ Hence
\begin{equation}\label{aux2}
\dot \Phi(t)^TJ_0 = \Phi(t)^T S(t).
\end{equation}
Using the definition of $S$ we compute
$$
S(t)^T =  (\Phi(t)^T)^{-1}\dot \Phi(t)^T J_0 \Rightarrow \dot \Phi(t)^T J_0 = \Phi(t)^T S(t)^T.
$$ 

It follows from \eqref{aux2} that  $S(t)^T=S(t), \forall t$. The reader can easily check that \eqref{periodicity} implies $S(t+1)=S(t) \forall t.$  Finally, \eqref{operator_unitary} follows from \eqref{Liederivative}, \eqref{definitionS} and the identity $$\mathcal{L}_S = -J_0 \cdot \Psi \circ \mathcal{L}_{\dot x_T} \circ \Psi^{-1}.$$

We check the symmetry of $A_{P,J}$ using coordinates $\Psi$.  The $L^2$-product on $\Gamma(x_T^*\xi)$ induced by $d\lambda$ and $J$ takes the form
$$
\left< \zeta_1,\zeta_2\right > = \int_0^1 \left<\zeta_1(t),\zeta_2(t) \right>\, dt, \ \ \ \ \forall \zeta_1,\zeta_2\in L^2(\R / \Z,\R^2).
$$
Integrating by parts, we obtain 
$$
\begin{aligned}
\left< \zeta_1,\mathcal{L}_S\cdot  \zeta_2\right> & = \int_0^1 \left<\zeta_1(t),-J_0\dot \zeta_2(t) - S(t)\zeta_2(t)\right>\, dt\\
& = \int_0^1 \left<-J_0\dot \zeta_1(t) - S(t)\zeta_1(t),\zeta_2(t)\right>\, dt\\
& = \left<\mathcal{L}_S \cdot \zeta_1,\zeta_2\right>.
\end{aligned}
$$ 

The eigenvalues of $\mathcal{L}_S$ are real.  A non-trivial $\nu$-eigenfunction $\zeta_\nu$ of $\mathcal{L}_S$  satisfies the smooth linear ODE
$$
\dot \zeta_\nu(t)= J_0 (S(t)+\nu I) \zeta_\nu(t) , \ \ \ \forall t,
$$  
and thus $\zeta_\nu$ is smooth and never vanishes. In particular, $\zeta_\nu$ has a well-defined winding number
$$
 \wind (\nu) = \frac{\Theta(1) - \Theta(0)}{2\pi}
$$ 
where  $\zeta_\nu(t) = (r(t)\cos(\Theta(t)),r(t) \sin(\Theta(t))) \forall t,$ for continuous functions $r>0,\Theta$. It does not depend on the $\nu$-eigenfunction.

If the linearized first return map $D\varphi_T(x_T(0)):\xi|_{x_T(0)} \to \xi_{x_T(0)}$  is the identity map then it is always possible to choose $J\in \mathcal{J}_+(\xi)$ and a trivialization $\Psi$ so that $S \equiv 0$. In this case, the eigenvalues of $\mathcal{L}_0=-J_0\frac{d}{dt}$ are
$$
\nu_k = 2\pi k, k \in \Z.
$$ 
Each $\nu_k$ admits a $2$-dimensional eigenspace
$$
\{t\mapsto A(\cos(2 \pi kt),\sin(2\pi kt)) + B(\sin(2\pi kt),-\cos(2\pi kt)), t\in \R / \Z,A,B\in \R\},
$$ 
whose eigenfunctions have winding number $k$.

The next theorem asserts that the general case is similar.
\begin{thm}[\cite{props2}]Let $t\mapsto S(t)$ be a smooth loop of $2\times 2$ symmetric matrices and let $\mathcal{L}_S$ be defined as in \eqref{operator_unitary}. Then
\begin{itemize}
\item[(i)] the spectrum $\sigma(\mathcal{L}_S)$ consists of real eigenvalues which accumulate precisely at $\pm \infty$.
\item[(ii)] the winding number $\wind(\nu)\in \Z,\nu\in \sigma(\mathcal{L}_S),$ is independent of the $\nu$-eigenfunction.
\item[(iii)] the map 
$$
\nu \mapsto \wind(\nu)\in \Z,
$$ is a surjective increasing map. For every $k\in \Z$, there exist precisely two eigenvalues $\nu_k^1,\nu_k^2$, counting multiplicities, so that $$\wind(\nu_k^1)=\wind(\nu_k^2)=k.$$ 
\item[(iv)] if $\nu_1\neq \nu_2$ satisfy $\wind(\nu_1)=\wind(\nu_2)$, then any two  $\nu_1,\nu_2$-eigenfunctions  are pointwise linearly independent.
\end{itemize}
\end{thm}

\section{The generalized Conley-Zehnder index}
\label{appendix_CZrotation}

Let $P=(x,T)$ be a closed Reeb orbit of $\lambda$ and let $\tau$ be a unitary trivialization of $x_T^*\xi$. Let $\mathcal{L}_S$ be the operator defined in \eqref{operator_unitary}. Let
$$
\begin{aligned}
\nu^\tau_-(P):= \max \{ \nu:  \nu \in \sigma(\mathcal{L}_S)\cap (-\infty,0)\},\\
\nu^\tau_+(P):= \min \{ \nu:  \nu \in \sigma(\mathcal{L}_S)\cap[0,+\infty)\},
\end{aligned}
$$
and let
$$
\wind_+^\tau(P):= \wind(\nu^\tau_+(P)) \ \ \ \mbox{ and } \ \ \  \wind^\tau_-(P):= \wind(\nu^\tau_-(P)) .
$$
These winding numbers satisfy $0\leq \wind^\tau_+(P)-\wind^\tau_-(P)\leq 1$ and they do not depend on  $J$. 

\begin{defn}[\cite{convex}] The generalized Conley-Zehnder index of $P=(x,T)$ with respect to  $\tau$ is defined as
$$
\cz^\tau(P):=\wind_+^\tau(P)+\wind_-^\tau(P).
$$
\end{defn}

This definition immediately implies that if $\nu\in \sigma(A_{P,J})$, then
$$
\begin{aligned}
\nu  < 0 \ \Rightarrow \wind^\tau (\nu) \leq \frac{\cz^\tau(P)}{2},\\ 
\nu \geq  0 \ \Rightarrow \wind^\tau (\nu) \geq \frac{\cz^\tau(P)}{2}.
\end{aligned}
$$

A more geometric definition of $\cz^\tau$ is as follows. Consider $\Phi_t,t\in[0,1]$, the family of $2\times 2$ symplectic matrices representing the linearized Reeb flow on $\xi$ along $P$ in coordinates induced by $\tau$ as in  the previous section. 

For each initial condition $0\neq \bar \zeta$ let $\Theta_{\bar \zeta}(t)$ be a continuous argument of $\Phi_t\bar \zeta$ with $t\in[0,1]$. Let
$$
\Delta\Theta (\bar \zeta) = \frac{\Theta_{\bar \zeta}(1) - \Theta_{\bar \zeta}(0)}{2\pi},
$$ 
and
$$
I^\tau(P):= \{ \Delta \Theta(\bar \zeta):0\neq \bar \zeta \in \R^2\},
$$ 
be the interval containing the argument variations of all initial conditions.

The length of $I^\tau(P)$ is less than $\frac{1}{2}$ and thus, for each $\epsilon>0$ small, either $I^\tau(P)-\epsilon$ contains an integer $k$ or is contained in between two consecutive integers $k$ and $k+1$. In the first case we define $\widetilde \mu^\tau(P) := 2k$ and in the second case, $\widetilde \mu^\tau(P)= 2k+1$.

It is a simple exercise to check that
$$
\cz^\tau(P) = \widetilde \mu^\tau(P).
$$
 Moreover, $P$ is nondegenerate if and only if the boundary of $I^\tau(P)$ does not contain an integer. For a proof to these facts, see \cite{fols}.

\section{The quaternionic trivialization}
\label{ap_quaternios}

Let $S \subset \R^4$ be a regular energy level of a Hamiltonian function $H$, that is $S=H^{-1}(c), c \in \R,$ and $\nabla H|_S$ never vanishes. The quaternion group induces an orthonormal frame
\begin{equation}\label{eq_frame}
TS = {\rm span}\{X_1,X_2,X_3\}|_S,
\end{equation}
spanned by the vector fields
\begin{equation}\label{eq_Xi}
X_i=A_i \dfrac{\nabla H}{|\nabla H|} \subset TS, \,\,\, i=1,2,3.
\end{equation}
 Here,  the $4\times 4$ matrices
$A_i,i=1,2,3,$ are
$$
A_1 = \left(\begin{array}{cc}  0 & J \\ J & 0 \end{array} \right) \ \ \  A_2 =\left(\begin{array}{cc} J & 0 \\ 0 & -J \end{array} \right) \ \ \ A_3 = \left(\begin{array}{cc} 0 & I \\ -I & 0 \end{array} \right),
$$
where $0$, $I$ and $J$ are the $2 \times 2$ matrices
\begin{equation}\label{eq_J}
0= \left(\begin{array}{cc} 0 & 0 \\ 0 & 0 \end{array} \right) \ \ \ \ I = \left(\begin{array}{cc} 1 & 0 \\ 0 & 1 \end{array} \right)  \ \ \ \ J = \left(\begin{array}{cc} 0 & 1 \\ -1 & 0 \end{array} \right).
\end{equation}
Observe that $X_3$ is parallel to the Hamiltonian vector field $X_H=A_3 \nabla H$.

Denote by $\phi_t$ the Hamiltonian flow of $X_H$ restricted to $S$. Since $d\phi_t: TS \to TS$ preserves the line bundle $\R X_3$, one may restrict the study of the linearized flow to 
$$TS/\R X_3 \simeq {\rm span} \{X_1,X_2\}.$$

As discussed in  Section \ref{sec_basics}, if $S$ has contact-type then $X_H|_S$ is parallel to a Reeb vector field $R$ of a contact form $\lambda$ on $S$. In this case, the contact structure $\xi=\ker\lambda$ is transverse to $R$ and as a result we have $\xi \simeq {\rm span} \{X_1,X_2\}$. This makes the quaternionic trivialization a useful tool for estimating Conley-Zehnder indices.

The orthonormal frame \eqref{eq_frame} induces a trivialization $\Psi: TS \to S \times \R^3$
\begin{equation}\label{eq_triv}
\Psi: TS \ni a_1X_1 + a_2X_2 +a_3X_3 \mapsto (a_1,a_2,a_3)\in \R^3
\end{equation}
which provides a simple form to analyze the transverse linearized flow along Hamiltonian trajectories.

\begin{prop}\label{prop_lin} In coordinates $(a_1,a_2)\in \R^2$ induced by the trivialization \eqref{eq_triv}, a solution to the linearized flow along a non-constant trajectory of $X_H$, projected to $TS/\R X_3$, satisfies the equation
\begin{equation}\label{linear}
\left(\begin{array}{c}\dot  a_1 \\ \dot a_2 \\  \end{array} \right) =-J M\left(\begin{array}{c} a_1 \\ a_2 \end{array} \right),
\end{equation}
where $J$ is given in \eqref{eq_J}, $M$ is the symmetric matrix
$$
M = \left(\begin{array}{cc} \kappa_{11} + \kappa_{33} & \kappa_{12} \\ \kappa_{12} & \kappa_{22}+\kappa_{33} \end{array} \right),
$$
and
$$
\kappa_{ij} = \left< {\bf H}  X_i,X_j \right>.
$$
The matrix ${\bf H}={\bf H}(x)$ is the Hessian of $H$ at $x\in S$, $X_i, i=1,2,3$, is given by \eqref{eq_Xi}, and $\left<\cdot,\cdot\right>$ is the standard inner product on $\R^4$.
\end{prop}

\proof 
Let $x(t)\in S$ be a Hamiltonian trajectory, that is a solution to 
\begin{equation}\label{eq_ham}
\dot x = A_3 \nabla H(x).
\end{equation}
A solution $y(t)\in T_{x(t)}S$ to the linearized flow $d\phi_t:TS\to TS$  along $x(t)$ satisfies the linear differential equation
\begin{equation}\label{linear2}
\dot y = A_3 {\bf H}(x) y.
\end{equation} Substituting $y=a_1X_1 + a_2X_2+a_3X_3$ in \eqref{linear2}, we obtain
\begin{equation}\label{eq_sum}
\sum_{i=1}^3 (\dot a_iX_i+a_i\dot X_i) = \sum_{i=1}^3 a_i A_3 {\bf H} X_i.
\end{equation}

We may assume for simplicity that $|\nabla H|=1$. In particular, it follows from \eqref{eq_Xi} and \eqref{eq_ham} that
\begin{equation}\label{Nponto}
\dot X_i = A_i {\bf H} \dot{x}=A_i {\bf H}A_3 \nabla H=A_i{\bf H}X_3.
\end{equation}
Taking the inner product of the expression \eqref{eq_sum} with $X_1$, using \eqref{Nponto} and the relations
$$
\begin{aligned}
A_i^T &=-A_i , \ \ \  
A_1A_2 = A_3, \ \ \  
A_2A_3 = A_1, \ \ \  A_3A_1=A_2, \ \ \  
\\
\left< X_i,X_j\right> &= \delta_{ij} \Rightarrow \left<\dot X_i, X_i\right>=0,  \ \ \  \left<\dot X_i, X_j\right> = - \left<\dot X_j,X_i\right> \ \ \ \forall i,j,
\end{aligned}
$$ we obtain
$$
\dot a_1 =- a_2\left<{\bf H}X_3,X_3\right>   -a_1\left< {\bf H}X_1,X_2\right> - a_2\left<{\bf H}X_2,X_2 \right>.
$$ This is precisely the expression for $\dot a_1$ in \eqref{linear}. Analogously one obtains $$
\dot a_2 =a_1\left<{\bf H}X_1,X_1\right>   +a_1\left< {\bf H}X_3,X_3\right> + a_2\left<{\bf H}X_1,X_2 \right>.
$$
\qed

\end{document}